\documentclass[12pt,a4paper]{amsart}

\input prepictex.tex
\input pictex.tex
\input postpictex.tex
\chardef\bslash=`\\ 

\makeatletter
\def\verbatim{\interlinepenalty\@M \@verbatim
  \leftskip\@totalleftmargin\advance\leftskip2pc
  \frenchspacing\@vobeyspaces \@xverbatim}
\makeatother
\newtheorem{theorem}{Theorem}[section]
\newtheorem{corollary}[theorem]{Corollary}
\newtheorem{lemma}[theorem]{Lemma}
\newtheorem{proposition}[theorem]{Proposition}

\newtheorem{remark}[theorem]{Remark}
\newtheorem{example}[theorem]{Example}
\theoremstyle{definition}
\newtheorem{definition}{Definition}[theorem]

\numberwithin{equation}{section}

\newcounter{picture}

\newcommand{\supp}{\operatorname{supp}}

\newcommand{\ZZ}{{\Bbb Z}}
\newcommand{\cA}{{\mathcal A}}
\newcommand{\D}{{\Delta}}
\newcommand{\G}{{\Gamma}}
\newcommand{\Om}{{\Omega}}

\newcommand{\w}{{\omega}}

\newsymbol\rtimes 226F 

\begin{document}

\title[MASA's of group factors]{ Maximal abelian subalgebras of the group
factor
of an $\widetilde A_2$ group }

\date{August 1, 1996}
\author[G. Robertson]{Guyan Robertson}
\author[T. Steger]{Tim Steger}
\address{Department of Mathematics  \\
        University of Newcastle\\  NSW  2308\\ AUSTRALIA}
\email{guyan@@maths.newcastle.edu.au }
\address{Istituto di Matematica e Fisica \\
Universit\`a di Sassari \\
via Vienna 2 \\ 07100 Sassari \\ ITALIA}
\email{steger@@ssmain.uniss.it }
\thanks{This research was supported by the Australian Research Council.} 
\thanks{ \hfill Typeset by  \AmS-\LaTeX}

\maketitle
\begin{abstract}
An $\widetilde A_2$ group $\Gamma$ acts simply transitively on the
vertices of an affine building $\triangle$. We study certain subgroups
$\Gamma_0 \cong {\Bbb Z}^2$ which act on certain apartments of
$\triangle$.  If one of these subgroups acts simply transitively on an
apartment, then the corresponding subalgebra of the group von Neumann
algebra is maximal abelian and singular. Moreover the Puk\'anszky
invariant contains a type $I_{\infty}$ summand.
\end{abstract}

\section*{Introduction}

Let $\Gamma$ be a group acting simply transitively on the vertices of a
homogeneous
tree $\bf T$ of degree $n+1<\infty$.  Then, by \cite[Chapter I, Theorem
6.3]{ftn},
$$
\Gamma\cong{\Bbb Z} * \cdots * {\Bbb Z} * {\Bbb Z}_2 * \cdots * {\Bbb Z}_2
$$
where there are $s$ factors of ${\Bbb Z}$, $t$ factors of ${\Bbb Z}_2$, and
$2s+t=n+1$.

We can identify $\bf T$ with the Cayley graph of $\Gamma$ constructed
from right multiplication by the natural generators for $\Gamma$. The
action of $\Gamma$ on $\bf T$ is equivalent to the natural action of
$\Gamma$ on its Cayley graph via left multiplication. With this
geometric interpretation, certain geodesics in $\bf T$ arise as Cayley
graphs of subgroups of $\Gamma$. That is, there are subgroups of
$\Gamma$ which act simply transitively on geodesics in $\bf T$.

\begin{example} \label{0}
{\em $\Gamma\cong{\Bbb Z}*{\Bbb Z}_2* {\Bbb Z}_2$ has generators
$a$,$b$,$c$, and relations $b^2 = c^2 = e$. Here $\bf T$ is a
homogeneous tree of degree four.  The Cayley graph of the subgroup
$\Gamma_0 = \langle a \rangle \cong {\Bbb Z}$ relative to the
generators $a, a^{-1}$ is a geodesic in $\bf T$.  The Cayley graph of
the nonabelian subgroup $\Gamma_1 = \langle b,c \rangle \cong{\Bbb
Z}_2 *{\Bbb Z}_2$ relative to the generators $b,c$ is also a geodesic.
}
\end{example}

We denote by $W^*(\Gamma)$ the group von Neumann algebra of a group
$\Gamma$.
It is generated by the left regular representation of $\Gamma$ on
$l^2(\Gamma)$.
We regard $W^*(\Gamma)$ as a convolution algebra embedded in $l^2(\Gamma)$
\cite[Vol.II, 6.7]{kr}.
Then $W^*(\Gamma) = \{ f \in  l^2(\Gamma) : f*l^2(\Gamma)  \subset
l^2(\Gamma) \} $.
If $\Gamma_0$ is a subgroup of $\Gamma$ then $W^*(\Gamma_0)$ may be
identified 
with the set of functions in $W^*(\Gamma)$ whose support is contained in 
$\Gamma_0$.

In the example above, $W^*(\Gamma_0)$ is an abelian von Neumann subalgebra
of $W^*(\Gamma)$. It is actually a maximal abelian subalgebra or {\it masa}.
In fact it follows from \cite [Proposition 4.1]{po1} that it is a {\it
singular} masa.
(The definition will be given later.)

The homogeneous tree $\bf T$ may be regarded as a one dimensional affine
building of type 
$\widetilde A_1$. Our purpose is to generalize the above observations
to the two dimensional case where $\Gamma$ is a group acting simply
transitively
 on the vertices of an affine building $\triangle$ of type $\widetilde
A_2$.
The building $\triangle$ is a chamber system consisting of vertices, edges
and triangles. 
 Each edge lies on $q + 1$ triangles, where $q \geq 2$ is the {\it order}
of
$\triangle$.  An {\it apartment} is a subcomplex of $\triangle$ isomorphic
to the Euclidean
plane tesselated by equilateral triangles (i.e. a Coxeter complex of type
$\widetilde A_2$).

An $\widetilde A_2$ group 
can be constructed as follows \cite [I,\S3]{cmsz}.
Let $(P,L)$ be a projective plane of order $q$. There are $q^2 + q + 1$
points
 (elements of $P$) and $q^2+q+1$
lines (elements of $L$).  Each point lies on $q+1$ lines and each line
contains $q+1$ points. 
Let $\lambda : P \rightarrow L$ be a bijection (a {\it point--line
correspondence}).  Let ${\mathcal T}$
be a set of triples $(x, y, z)$ where $x, y, z \in P$, with the following
properties.

(i)  Given $x, y \in P$, then $(x, y, z) \in {\mathcal T}$ for some
$z \in P$ if and only if $y$ and $\lambda(x)$ are incident (i.e. $y \in 
\lambda(x)$).

(ii)  $(x, y, z) \in {\mathcal T} \Rightarrow (y, z, x) \in {\mathcal T}$.

(iii)  Given $x, y \in P$, then $(x, y, z) \in {\mathcal T}$ for at
most one $z \in P$.

${\mathcal T}$ is called a {\it triangle presentation}  compatible with
$\lambda$.  A complete list is given in \cite{cmsz} of all triangle
presentations for $q = 2$
and $q = 3$.

Let $\{a_x : x \in P\}$ be  $q^2 + q + 1$ distinct letters and form the
group
$$
\Gamma = \big\langle a_x, x \in P \ |\  a_x a_y a_z = e \hbox { for } (x,
y, z) \in {\mathcal T}
\big \rangle $$

 The Cayley  graph of $\Gamma$ with respect to the
generators $a_x, x \in P$, and their inverses is the $1$-skeleton of an 
affine building of type $\widetilde A_2$.

It is convenient to identify the point $x \in P$ with the generator
$a_x \in \Gamma$.  The lines in $L$ correspond to the inverses of the
generators of $\Gamma$ according to the point-line correspondence:
$a_{\lambda(x)} = a_x^{-1}$ for $x \in P$ \cite{cmsz}. We may
therefore write $x^{-1}$ for $a_x^{-1}$ and identify $x^{-1}$ with
$\lambda(x)$. From now on we use the concise notation $x$ and
$\lambda(x)$ instead of $a_x$ and $a_{\lambda (x)}$ respectively. It
is important to note that ${\mathcal T}$ can be recovered from $\Gamma$ :
$${\mathcal T} = \{ (x,y,z) : x,y,z \in P  \hbox { and }   \ xyz = e \}$$
This implies that if $x,y \in P$ then $y \in \lambda(x)$ if and only if 
$xyz = e$ for some $z \in P$.

 Any element
$g \in \Gamma \setminus \{ e \}$ can be written uniquely in the {\it left
normal form}
$$g = x_1^{-1}x_2^{-1} \dots x_m^{-1}y_1y_2 \dots y_n$$
where there are no obvious cancellations and  $x_i,y_j \in P$, $1 \le i \le
m, 1 \le j \le n$
 \cite [Lemma 6.2]{cm}.
The absence of ``obvious'' cancellations means that
$ x_i \notin \lambda(x_{i+1})  \quad ( 1 \le i < m)$ ,
$ y_{j+1} \notin \lambda(y_j)  \quad (1 \le j < n)$, and 
$ x_m \ne y_1 $.
Also any such word for $g$ is a minimal word for $g$ in the generators $x
\in P$ and their inverses 
\cite [Lemma 6.2]{cm}. We write $|g| = m + n$.
 An exactly analogous statement is true for the {\it right normal
form} for $g$ in which the  inverse generators are on the right of the word
\cite 
[I,Proposition 3.2]{cmsz}. We shall use these facts repeatedly. The reader
is referred to \cite{kr} for background information on von Neumann algebras
and to \cite{bro}, \cite{ron} for buildings. Operator algebras associated
with $\widetilde A_2$
buildings are studied in \cite{ros} and \cite{rr}. S. Mozes \cite{moz} has
also been
concerned with automorphism groups of affine buildings and corresponding
actions
of ${\Bbb Z}^2$ on apartments. 

 From now on, unless otherwise stated, $\Gamma$ will denote an $\widetilde
A_2$
 group with associated projective plane of order $q \ge 2$, and $\triangle$
will denote
the corresponding affine building whose vertices are identified with the 
elements of $\Gamma$. The following result shows
that $W^*(\Gamma)$ is a factor \cite [Theorem 6.7.5]{kr}.

\begin{lemma} \label{1} $\Gamma$ is an i.c.c. group. That is, each
conjugacy class
in $\Gamma$, except for $\{ e \}$, is infinite.
\end{lemma}

{\sc Proof:} Let $g \in \Gamma \setminus \{e\}$. Assume that $g$ has left
normal form
 $$g = x_1^{-1}x_2^{-1} \dots x_m^{-1}y_1y_2 \dots y_n$$
where $m,n \ge 1$. ( If $m = 0$ or $n = 0$ the argument is simpler.)
Thus $|g| = m+n$.
We may choose $z\in P$ such that $z \notin \lambda(x_1)$ and $z \notin
\lambda(y_n)$.
This is possible because there are $q+1$ points on each of the lines
$\lambda(x_1)$ and $\lambda(y_n)$.
Since there are $q^2 + q + 1$ points altogether,
and $q^2 + q + 1 > 2(q+1)$, we may choose a point $z$ not lying on
either of these lines.
Then $z^{-1}x_1^{-1}x_2^{-1} \dots x_n^{-1}y_1y_2 \dots y_m z$ is in left
normal form.
It follows that $|z^{-1}gz| = |g| + 2$.
Repeating the process, we see that the conjugacy class of $g$
contains a sequence of elements of length $|g| + 2n, n = 1,2,\dots$
and hence is infinite.
\qed

\bigskip
 The $\widetilde A_2$ groups have Kazhdan's property (T) \cite [Theorem
4.6]{cms}.
It is therefore particularly interesting to investigate how properties
of such groups $\Gamma$ are reflected in the structure of $W^*(\Gamma)$,
in view of the following rigidity conjecture of A. Connes
\cite[V.B.$\varepsilon$]{con}.

{\it Conjecture: If $\Gamma_1$ and $\Gamma_2$ are i.c.c. groups 
with property (T) and $\Gamma_1$ is not isomorphic to $\Gamma_2$,
 then $W^*(\Gamma_1)$ is not isomorphic to $W^*(\Gamma_2)$.}
\bigskip

\section{Some abelian subgroups of $\widetilde A_2$ groups}

 Recall that an  apartment is a subcomplex of $\triangle$ isomorphic to the
Euclidean
plane tesselated by equilateral triangles. An abelian subgroup of $\Gamma$
which acts simply transitively on an apartment will be the analogue of the
subgroup
${\Bbb Z}$ in the tree case (Example \ref{0}).
In the $\widetilde A_2$ case, such an abelian subgroup $\Gamma_0$
necessarily
contains three distinct elements $a,b,c$ of $P$. We begin by elucidating
the structure of $\Gamma_0$.

\begin{lemma} \label{2} If $x,y \in P$, $x \ne y$ and $xy =yx$
then $xyz =e$, where $\{z\} = \lambda(x) \cap \lambda(y)$. Moreover
 $x \in \lambda(y)$ and $y \in \lambda(x)$.
\end{lemma}
{\sc Proof:} Suppose that $xy = yx$. Then the left side of the equation
$y^{-1}xy = x$ is not in left normal form, by the uniqueness assertion of
\cite[Lemma 6.2]{cm}. Since $y \ne x$ it follows that $y \in \lambda(x)$.
Thus $xyz = e$ for some $z \in P$. The fact that $yxz = e$ then shows that
$x \in \lambda(y)$. It is also clear from these equations that
$z \in \lambda(x) \cap \lambda(y)$ 
\qed

\begin{lemma} \label{a,b,c}
Let $a,b,c$ be distinct mutually commuting  generators in $P$. Then $abc =
e$
and $\lambda(a) \cap \lambda(b) = \{c\}$.
\end{lemma}
{\sc Proof:} By Lemma \ref{2} we have $abz =e$, where $\{z\} = \lambda(a)
\cap \lambda(b)$.
Since $ac = ca$, the same lemma shows that $c \in \lambda(a)$.
Since $bc = cb$, we also have $c \in \lambda(b)$.
Therefore $z = c$, which proves the result.
\qed  

\begin{remark}
{\em It follows from this result that a set of pairwise commuting
elements in~$P$ can have no more than three elements.  }
\end{remark}

\begin{lemma} \label{z2}
Let $a,b,c$ be distinct mutually commuting
 generators in $P$.
Then the subgroup $\Gamma_0 = \langle a,b,c \rangle$ of $\Gamma$  is
 isomorphic to ${\Bbb Z}^2$.
The Cayley graph of $\Gamma_0$ relative
to the generators $a,b,c$ and their inverses is the 1-skeleton of an
apartment in 
$\triangle$.
\end{lemma}

{\sc Proof:} By Lemma \ref{a,b,c}, $abc = e$.
  Therefore an isomorphism $\theta : \langle a,b \rangle \to {\Bbb Z}^2$ is
defined by
$\theta(a) = (1,0)$, $\theta(b) = (0,1)$, $\theta(c) = (-1,-1)$. The
vertices adjacent to $e$
in the corresponding Cayley graph are labelled as in Figure 1.
\qed

\refstepcounter{picture}
\begin{figure}[htbp]
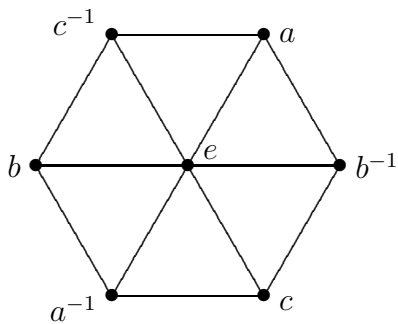
\label{fig1}
{}\hfill
\beginpicture
\setcoordinatesystem units <1cm, 1.732cm>  
\setplotarea  x from -2.5 to 2.5,  y from -1.5 to 1.5
\put {$\bullet$} at 0 0
\put {$\bullet$} at 1 1
\put {$\bullet$} at 1 -1
\put {$\bullet$} at -1 1
\put {$\bullet$} at -1 -1
\put {$\bullet$} at 2 0
\put {$\bullet$} at -2 0
\put {$e$} [b,l] at 0.2 0.05
\put {$a$} [l] at  1.2 1
\put {$b$} [r] at -2.2 0
\put {$c$} [t,l] at  1.2 -1
\put {$a^{-1}$} [r,t] at  -1.2 -1
\put {$b^{-1}$} [l] at  2.2 0
\put {$c^{-1}$} [b,r] at  -1.2  1
\putrule from -1 1 to 1 1
\putrule from -2 0 to 2 0
\putrule from -1 -1 to 1 -1
\setlinear \plot 0 0  1 1  2 0  1 -1  0 0  -1 1  -2 0  -1 -1  0 0 /
\endpicture
\hfill{}
\caption{The star of $e$ in the Cayley graph of ${\Bbb Z}^2$.}
\end{figure}

\bigskip
\begin{remark} \label{abelian subgroups}
{\em
There are many examples of such subgroups. For example, the groups (4.1),
(5.1), (6.1), (9.2), 
(13.1), (28.1) in the tables at the end of \cite{cmsz} contain such
subgroups. 
(However there are no examples when $q = 2$.)

 Note that each element  $g$ of the group $\langle a,b,c \rangle$ 
has left (right) normal form
$g = x^{-k}y^l \ (= y^lx^{-k})$ where $x,y \in \{a,b,c\}$, $k,l \ge 0$ and
$x \ne y$.
}
\end{remark}

\begin{lemma}\label{z2a}
Let $a,b \in P$ with $ab^2 =e$. If $a \ne b$ then  $\langle a,b \rangle
\cong {\Bbb Z}$
and  Cayley graph of the group $\langle a,b \rangle$ relative
to the generators $a,b$ and their inverses is an infinite strip in an
apartment of $\triangle$.
If $a = b$ then $\langle a,b \rangle = \langle a \rangle$ is cyclic of
order 
three and the Cayley graph is a triangle. 
\end{lemma}

{\sc Proof:}  
If $a \ne b$, an isomorphism $\theta : \langle a,b \rangle \to {\Bbb Z}$ is
defined by
$\theta(a) = 2$, $\theta(b) = -1$. It is easy to verify the remaining
assertions.
\qed

\begin{remark}
{\em
Lemma \ref{z2a} is a degenerate analogue of the Lemma \ref{z2}.
The case when $\Gamma_0 \cong {\Bbb Z}$ and the Cayley graph is an 
infinite strip occurs for subgroups of many of the groups given explicitly
in \cite{cmsz}.
 For example, the groups (B3)  (where $q = 2$)  and 
(9.3),(38.1),(63.1),(64.1) (where $q = 3$). A counting argument
\cite[II, \S3]{cmsz}
shows that if $q$ is divisible by $3$ then $\Gamma$ has a generator of
order three.  
}
\end{remark}

\begin{remark} \label{nonabelian}
{\em
In the case of a group $\Gamma$ acting simply transitively on the vertices of a tree
there may exist a nonabelian subgroup of $\Gamma$ 
(necessarily isomorphic to ${\Bbb Z}_2 *{\Bbb Z}_2$)
which acts simply transitively on a geodesic in the tree. 
 In Example \ref{0} the Cayley graph of the subgroup
 $\langle b,c \rangle \cong{\Bbb Z}_2 *{\Bbb Z}_2$
relative to the generators $b,c$ and their inverses is a geodesic in the
tree.

 Analogously, a nonabelian subgroup of an $\widetilde A_2$ group 
can act simply transitively on an apartment.
There are two examples.
  One is the Dyck group (or triangle group) 
$$T(3,3,3) = \langle x, y, z\  |\  x^3 = e, y^3 = e, z^3 = e, xyz = e
\rangle$$
$T(3,3,3)$ is a subgroup of index $2$ of the Coxeter group $W$ of type
$\widetilde A_2$.
 In fact $T(3,3,3)$ is the
``rotation subgroup'' of $W$ consisting or words of even length in the
canonical generators of $W$ \cite[II.3 and II.4]{mag}.

 The other possibility is the amalgam
$${\Bbb Z} *_{\Bbb Z}{\Bbb Z} = \langle  y, z\  |\  y^2 = z^2 \rangle
= \langle x, y, z\  |\  xy^2 = xz^2 = e \rangle$$
This is the nonabelian Bieberbach group in two dimensions, namely
the fundamental group of the Klein bottle \cite[Chapter 1]{ch}.
The Cayley graph of each of the groups $T(3,3,3)$ and ${\Bbb Z} *_{\Bbb
Z}{\Bbb Z}$
relative to the generators $x$,$y$,$z$ and their inverses  is
the 1-skeleton of the Coxeter complex of type $\widetilde A_2$.

Many of the $\widetilde A_2$ groups enumerated in \cite{cmsz} contain 
triples of
generators which generate a subgroup isomorphic to one of these groups.
Detailed enumeration of the possibilities shows that the groups ${\Bbb
Z}^2$,
presented as in Lemma \ref{z2}, and  the  groups
$T(3,3,3)$, ${\Bbb Z} *_{\Bbb Z}{\Bbb Z}$, presented as 
above, are the only possible subgroups of an
$\widetilde A_2$ group which can act simply transitively on an apartment. 

More precisely, let $\G$ be a group of automorphisms of a Coxeter complex
of type $\widetilde A_2$ which acts simply transitively and in a type
rotating way on the vertices. Then the generators satisfy relations of the
form $xyz = e$ and the $1$-skeleton of the Coxeter complex is the
Cayley graph of $\G$.There are three generators $a$,$b$,$c$ and the
neighbours of $e$ lie as shown in Figure \ref{fig2}, where
$\{g_1,g_2,g_3\} = \{a^{-1},b^{-1},c^{-1}\}$.
\bigskip

\refstepcounter{picture}
\begin{figure}[htbp]
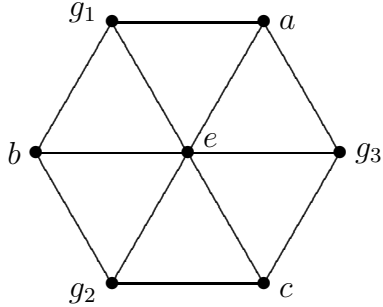
\label{fig2}
{}\hfill
\beginpicture
\setcoordinatesystem units <1cm, 1.732cm>  
\setplotarea  x from -2.5 to 2.5,  y from -1.5 to 1.5
\put {$\bullet$} at 0 0
\put {$\bullet$} at 1 1
\put {$\bullet$} at 1 -1
\put {$\bullet$} at -1 1
\put {$\bullet$} at -1 -1
\put {$\bullet$} at 2 0
\put {$\bullet$} at -2 0
\put {$e$} [b,l] at 0.2 0.05
\put {$a$} [l] at  1.2 1
\put {$b$} [r] at -2.2 0
\put {$c$} [t,l] at  1.2 -1
\put {$g_2$} [r,t] at  -1.2 -1
\put {$g_3$} [l] at  2.2 0
\put {$g_1$} [b,r] at  -1.2  1
\putrule from -1 1 to 1 1
\putrule from -2 0 to 2 0
\putrule from -1 -1 to 1 -1
\setlinear \plot 0 0  1 1  2 0  1 -1  0 0  -1 1  -2 0  -1 -1  0 0 /
\endpicture
\hfill{}
\caption{The star of $e$ in a Coxeter complex.}
\end{figure}
\bigskip

There is a relation of the form $xyz = e$ if and only if $yz = x^{-1}$,
that is $y$ is adjacent to $x^{-1}$ in the graph.
We now have a degenerate version of \cite[I, \S3]{cmsz}.
Namely, the link of $e$ is the incidence graph of a {\it geometry}
with set of points $P = \{a,b,c\}$ and set of lines 
$L = \{a^{-1},b^{-1},c^{-1}\}$. Moreover, $y \in x^{-1}
\Longleftrightarrow xyz =e$, for some generator $z$.
There are essentially three cases :
$(g_1,g_2,g_3)$ is one of the triples
$(c^{-1},a^{-1},b^{-1})$, $(a^{-1},b^{-1},c^{-1})$,
$(b^{-1},a^{-1},c^{-1})$.
Detailed enumeration of the possibilities shows that each geometry
gives rise to exactly one group whose Cayley graph tesselates 
the plane by equilateral triangles. The groups are respectively
${\Bbb Z}^2$, $T(3,3,3)$, and ${\Bbb Z} *_{\Bbb Z}{\Bbb Z}$.
It is instructive to sketch the labelled Coxeter complex in
each case. 
These groups may be regarded as degenerate 
$\widetilde A_2$ groups corresponding to a
degenerate projective plane $(P,L)$
of order $q = 1$.
}
\end{remark}

\bigskip
We now consider a much more general situation.  The apartment of Lemma
\ref{z2} which is spanned by two commuting generators of $\G$ is an
example of a periodic apartment.  An apartment $\cA$ in $\D$ is {\it
doubly periodic} if, when regarded as a copy of the Euclidean plane
tesselated by equilateral triangles, with directed edges labelled by
generators of $\G$, it has two independent periodic directions.  This
means that the edge-labelled apartment is invariant under an action of
$\ZZ^2$ on $\cA$ by translation.  This definition depends on the
choice of the group $\G$ acting on $\triangle$.

This concept  coincides with that of a {\it rigidly periodic } apartment in
the sense of \cite{rr}. (See \cite[Lemma 2.10]{rr}.)
A {\it rigidly periodic limit point} $\w \in \Om$ is a boundary
point of a rigidly periodic apartment.
It is important to note the elementary fact that each doubly periodic
apartment
 is uniquely determined by any one of its sectors \cite[Lemma 2.4]{rr}.

Now let $\cA$ be a doubly periodic apartment, and {\it assume
that $\cA$ contains the vertex $e$}.
 By double periodicity of the apartment we can choose independent periodic
directions
in $\cA$. This means that there are vertices $u$, $v$
 such that the edge-labelling of the apartment is
identical when viewed from each vertex $u^m$, $v^n$, $m,n \in \ZZ$. Then
$u^m\cA = \cA$, $v^n\cA=\cA$ and left multiplication by $u^mv^n$ defines
a periodic action of $\ZZ^2$  on $\cA$
by translation.
The {\it period group} $\G_0$ of a periodic apartment $\cA$ containing $e$
is the  abelian subgroup of $\G$ consisting of all vertices of the
apartment
which are periodic directions in the above sense. Clearly $\G_0$ is 
isomorphic to $\ZZ^2$.
By analogy with the tree case, one would expect $W^*(\Gamma_0)$ to be
a masa of $W^*(\Gamma)$.

\bigskip
\section{MASA's of $W^*(\Gamma)$}

Suppose that $\Gamma_0$ is a maximal abelian subgroup of a countable
discrete group $\Gamma$.  Suppose furthermore that the following two
conditions are satisfied.

(I) If $\G_0 x_0 \subset y_1\G _0 \amalg \dots \amalg y_n\G_0$,
where $y_1, \dots ,y_n \in \G$,
then $\G_0x_0 = y_j\G_0$ for some $j$.

(II) If $\phi$ is an automorphism of $\G_0$ which fixes pointwise some
finite index subgroup $K < \Gamma_0$, then $\phi$ is the identity
automorphism of $\G_0$.

\medskip
Conditions (I) and (II) are quite different in
nature: condition (I) depends on  how $\G_0$ lies in $\G$,
while condition (II) is a purely algebraic condition on $\G_0$. 
\bigskip
\begin{proposition}\label{masa} Under the above hypotheses,
  $W^*(\Gamma_0)$ is a masa of $W^*(\Gamma)$.
\end{proposition}
 
{\sc Proof:} Let $u \in W^*(\Gamma_0)'$ be unitary.
Suppose that $x_0 \in \supp u$. We show that $x_0 \in \G_0$.
This will prove that $u \in W^*(\Gamma_0)$.

For $g_0 \in \G_0$ we have $\delta_{g_0} \star u \star \delta_{g_0^{-1}} =
u$.

Therefore $u(g_0 x_0 g_0^{-1}) = u(x_0)$.
Since $u(x_0) \ne 0$ and $u \in l^2(\G)$, the set 
$\{ g_0 x_0 g_0^{-1} : g_0 \in \G_0 \}$ is finite.
Write its elements as $y_1, \dots, y_n$ and if necessary
delete repetitions which define the same cosets in $\G / \G_0$.
Then $\G_0x_0 \subset y_1\G_0 \amalg \dots \amalg y_n\G_0$.

By condition (I), $\G_0x_0 = y_j\G_0$ for some $j$.
Now $y_j = g_0 x_0 g_0^{-1}$ for some $g_0 \in \G_0$.
Thus $\G_0 x_0 = g_0x_0\G_0$. In other words,  
$\G_0 x_0 = x_0\G_0$. That is $x_0\G_0 x_0^{-1} = \G_0$.
We may therefore define an automorphism $\phi$ of
$\G_0$ by $\phi (g_0) = x_0g_0x_0^{-1}$.
Then for each $g_0 \in \G_0$  $\phi(g_0)x_0 = x_0g_0$
and so
$u(x_0) = u(g_0^{-1}x_0g_0) = u(g_0^{-1}\phi(g_0)x_0)$.

Again, since $u(x_0) \ne 0$ and $u \in l^2(\G)$, it follows that
the set $\{ g_0^{-1}\phi(g_0) : g_0 \in \G_0 \}$ is finite.
Let $K = \{ g_0 \in \G_0 : \phi (g_0) = g_0 \}$.
Then $K$ is the kernel of the homomorphism $g_0 \mapsto g_0^{-1}\phi(g_0)$
on $\G_0$, which has finite range. Therefore $\G_0/K$ is finite.
Thus $\phi$  fixes the finite index subgroup $K < \G_0$.

It follows from condition (II) that $\phi$ is the identity automorphism
of $\G_0$. This means that $x_0g_0 = g_0x_0$ for all $g_0 \in \G_0$.
Since we assumed that $\G_0$ is a maximal abelian subgroup of $\G$,
it follows that $x_0 \in \G_0$.
This proves that $u \in W^*(\Gamma_0)$ and hence
that $W^*(\Gamma_0)$ is a masa of $W^*(\Gamma)$.
\qed

\bigskip

\begin{lemma} \label{sectors}
Fix $C>0$ and let $S$ and~$S'$ be sectors (Weyl chambers)
in an affine building.  Either $S$ and~$S'$
share a common subsector, or $S$~has a subsector all of whose points
are at distance $>C$ from~$S'$.
\end{lemma}

{\sc Proof:} Choose subsectors $S_1$ and~$S_1'$ of $S$ and~$S'$
respectively which lie in a common apartment \cite[Chapter 9, Proposition (9.5)]{ron}.
If $S_1$ and~$S_1'$ point in the same direction, then they have a
common subsector, which is also a common subsector of~$S$ and~$S'$.

Otherwise, fix a finite $C_1>0$ so that $d(v,S_1')\leq C_1$ for
any~$v\in S'$ \cite[Chapter 9, Lemma (9.2)]{ron}.  Choose a subsector $S_2$
of~$S_1$ all of whose points are at distance $>C+C_1$ from $S_1'$.
Then those points are all at distance $>C$
from $S'$.
\qed
\bigskip

\begin{lemma} \label{extra} Let $\G$ be an $\widetilde A_2$ group 
and $\cA$ a doubly periodic apartment in the corresponding building.
Suppose that there exist $y_1, \dots ,y_n \in \G$ such that the
Hausdorff distance from $\cA$ to $y_1\cA \cup \dots \cup y_n\cA$ is
finite. Then $\cA$ coincides with some $y_j\cA$.
\end{lemma}


{\sc Proof:} Take a sector~$S$ in $\cA$.  Write each of the
apartments~$y_j\cA$ as a finite union of sectors~$S'$.  Suppose that
$S$~does not have a subsector in common with any of those
sectors~$S'$.  Fix $C>0$.  By Lemma~\ref{sectors}, some finite
intersection of subsectors of $S$ has all its points at distance $>C$
from $y_1\cA \cup \dots \cup y_n\cA$.  A finite intersection of
subsectors of~$S$ is nonempty (in fact, another subsector) so this
contradicts the main hypothesis.

Now we know that for some~$j$ the doubly periodic apartments $\cA$ and
$y_j\cA$ have a common subsector, and therefore that they coincide.
\qed
\bigskip

\begin{lemma} \label{lattice} Assume that $\G_0 \cong {\Bbb Z}^2$ is a
subgroup of the period group of a doubly periodic apartment $\cA$
containing the vertex $e$. Then conditions (I) and (II) hold.  If,
furthermore, $\G_0$ is the full period group of $\cA$, then $\G_0$ is
a maximal abelian subgroup of $\G$.
\end{lemma}

{\sc Proof:} (I) Suppose that
 $\G_0 x_0 \subset y_1\G _0 \amalg \dots \amalg y_n\G _0$,
 where $y_1, \dots ,y_n \in \G$. Then the Hausdorff distance
$d(\G _0,  y_1\G _0 \amalg \dots \amalg y_n\G _0) \le |x_0|$
and so there exists $C > 0$ such that
 $$d(\cA,  y_1\cA \amalg \dots \amalg y_n\cA) \le C$$
It follows from Lemma \ref{extra} that $\G_0x_0 = y_j\G_0$.

(II) Let $K$ be a subgroup of $\G_0 \cong {\Bbb Z}^2$. There is a basis 
$\{e_1, e_2\}$ of $\G _0$ such that $K = \langle  e_1^m |\  e_2^n \rangle$.
If $K$ has finite index, then $m \ne 0, n \ne 0$.
If $\phi$ is an automorphism of $\G_0$ which
fixes $K$ pointwise, then since $\Gamma_0$ is
torsion free,  $\phi(e_1) = e_1, \phi(e_2) = e_2$,
and so $\phi$ is the identity automorphism of $K$.

Now assume that $\G_0$ is
 the full period group  of $\cA$ and suppose that
an element $g \in \G$ commutes with $\G_0$.
We must show that $g \in \G_0$.
If $g_0 \in \G_0$ then $gg_0 = g_0g$ and so  $d(gg_0, g_0) \le |g|$.
Thus $$d(g\G _0, \G_0) \le |g|$$
Therefore there is a constant $C > 0$ such that
$$d(g\cA, \cA) \le C$$
It follows from Lemma \ref{extra} that $g\cA = \cA$.
The elements of $\G_0$  are vertices
of the apartment $\cA$ and since $g$ commutes with $\G _0$,
 $g$ acts on $\G_0$ by translation.
Hence $g$ must act on $\cA$ by translation, rather than
rotation or glide-reflection.
Thus $g$ is a period of the apartment $\cA$ and so
$g$ belongs to~$\G_0$.
\qed
\bigskip

\begin{corollary} \label{max}
Assume that $\G_0$ is the full period group of a doubly periodic 
apartment $\cA$ containing the vertex $e$.
Then $W^*(\G_0)$ is a masa of $W^*(\G)$.
\end{corollary}

\begin{remark} 
{\em
It was shown in \cite[\S 6]{cm} that the weak closure of the
algebra of biradial functions on $\G$ is a masa of $W^*(\G)$.
}
\end{remark}
\bigskip

\begin{definition} \label{regular} If $A$ is an abelian
 von Neumann subalgebra of a von Neumann
algebra $M$ then the {\it normalizer} $N(A)$  is the set of unitaries $u$
in $M$
such that $u^*Au = A$. The subalgebra $A$ is called {\it regular} 
(or a {\it Cartan subalgebra}) if $N(A)'' = M$.
It is called {\it singular} if $N(A)'' = A$.

Regular masa's are easy to construct from the classical group measure
space construction in which $M$ is the crossed product of an abelian
von Neumann algebra $A$ by a discrete group $G$ which acts ergodically
on $A$. Then $A$ is a regular masa in $M$. It is more difficult to
construct singular masa's. See \cite{rad} and \cite{br} for recent
examples.  The masa of biradial functions on an $\widetilde A_2$ group
has been studied in \cite{cms}.  The authors have convinced themselves
that this masa is singular.  Singularity has been proved in \cite{rad}
and \cite{br} for certain masa's of radial functions, in the context
of trees and similar rank~1 structures.

\end{definition}

\bigskip
\begin{lemma} \label{singular masa}
Suppose that $\Gamma$~is a maximal abelian subgroup of a countable
discrete group~$\Gamma$.  Assume in addition to conditions (I) and
(II) that no element of $\G \setminus \G_0$ normalizes $\G_0$.  Then
$W^*(\G_0)$ is a singular masa of $W^*(\G)$.
\end{lemma}

{\sc Proof:}  Let $u \in W^*(\Gamma)$ be a unitary satisfying
$u^*W^*(\G_0)u \subset W^*(\G_0)$. 
Suppose that $x_0 \in \supp u$. We must prove that $x_0 \in \G_0$.

There are only a finite number of cosets $y\G_0$ such that
$\|u|_{y\G_0}\|_2 \ge |u(x_0)|$.  Call them $y_1\G _0, \dots , y_n\G
_0$.  We claim that $\G_0 x_0 \subset y_1\G _0 \amalg \dots \amalg
y_n\G _0$.  To prove this, note that if $z \in \G_0$ then $u^{-1}
\star \delta_z \star u = f \in W^*(\G_0)$ is unitary. Therefore

\begin{multline}
|u(x_0)|
= |(\delta_z \star u)(zx_0)|
= |(u \star f)(zx_0)| \\
\le \|(u \star f)|_{zx_0\G_0}\|_2
=   \|u|_{zx_0\G_0} \star f\|_2
=   \|u|_{zx_0\G_0}\|_2
\notag
\end{multline}
(The last two equalities are valid because $\supp f \subset \G_0$
and  $f$ is unitary.)
This shows that $zx_0 \in y_1\G _0 \amalg \dots \amalg y_n\G _0$,
as claimed above.

Suppose that $x_0 \notin \G_0$. It follows from condition (I)
that $\G_0x_0 = y_j\G_0$ for some $j$.
In particular $x_0 \in y_j\G_0$, and so 
$x_0\G_0 = y_j\G_0 =\G_0x_0$.
Thus $x_0^{-1}\G_0x_0 = \G_0$ with  $x_0 \notin \G_0$,
contrary to the assumption of the lemma.  This contradiction shows
that $W^*(\Gamma_0)$ is singular.
\qed
\bigskip

\begin{theorem} \label{<a,b>singular masa}
Let $a,b,c$ be distinct mutually commuting generators of an
$\widetilde A_2$~group $\G$ and let $\Gamma_0 = \langle a,b,c \rangle
\cong {\Bbb Z}^2$.  (See Lemma \ref{z2}.)  Then $W^*(\G_0)$ is a
singular masa of $W^*(\Gamma)$.
\end{theorem}

{\sc Proof:} Note that the elements of $\G_0$  are the vertices
of a doubly periodic apartment $\cA$.
According to the preceding results,
we need only show that  no element of 
$\G \setminus \G_0$ normalizes $\G_0$.
Let $g \in \G$ satisfy $g\G_0g^{-1} = \G_0$.
If $g_0 \in \G_0$ then $gg_0 = g_0'g$ for some
$g_0' \in \G_0$ and so $d(gg_0, g_0') \le |g|$.
Thus $$d(g\G_0, \G_0) \le |g|$$
In other words,
$$d(g\cA, \cA) \le |g|$$
It follows from Lemma \ref{extra} that $g\cA = \cA$.
Since the elements of $\G_0$  are the vertices
of the apartment $\cA$, this means  $g\G_0 = \G_0$.
Therefore $g \in \G_0$, as required.
\qed

\begin{remark} \label{6} 
{\em
Suppose that $\G_0$ is the full period group of a doubly periodic 
apartment $\cA$ containing the vertex $e$, as in Corollary \ref{max}.
Then the masa $W^*(\G_0)$ of $W^*(\G)$ may be nonsingular.

For example in the group (B3) of \cite{cmsz}
there exist generators $a$, $b$, $c$ satisfying $ab^2 = ac^2 = e$ .
Let $\Gamma_0 = \langle a,bc \rangle \cong {\Bbb Z}^2$.
It is a subgroup of the group
$\langle a,b,c \rangle \cong {\Bbb Z} *_{\Bbb Z}{\Bbb Z}$
whose 1-skeleton is the Cayley graph of an apartment $\cA$.
(See Remark \ref{nonabelian}.)

Then $\Gamma_0$ is the period lattice of $\cA$
and so $W^*(\G_0)$ is a masa.
However $b \notin \G_0$ and $b$ normalizes $\G_0$. For 
$bab^{-1} = a \in \G_0$ and 
$b(bc)b^{-1} = c^3b^{-1} = cb^{2}b^{-1} = cb = a^{-1}(bc)^{-1}a^{-1} \in
\G_0$.
Therefore $W^*(\G_0)$ is not singular. 
}
\end{remark}

\begin{remark}
{\em
An $\widetilde A_2$ group $\Gamma$ has Kazhdan's property $T$ \cite{cms}.
It therefore follows from
\cite [Corollary 4.5]{po2} that $W^*(\Gamma)$ contains an {\it
ultrasingular} masa $A$. This
means that the only automorphisms of $W^*(\Gamma)$ which normalize  $A$ are
the inner
automorphisms implemented by unitaries in $A$. We do not know whether the
masa's we have
been considering are ultrasingular.
}
\end{remark}

\begin{remark}
{\em
H. Yoshizawa has given an explicit decomposition of the regular
representation
of the free group on two generators into two
inequivalent families of irreducible representations. 
These representations are induced from  characters on the abelian subgroup
generated by one of the generators. See, for example \cite[Chap.19]{rob}.
Our results may be re-interpreted to give a two dimensional analogue of
Yoshizawa's construction. For example, if $\G$, $\G_0$ are as in
Theorem \ref{<a,b>singular masa}, then $\widehat\Gamma_0 \cong {\Bbb T}^2$
and we have a direct integral decomposition
$$l^2(\G) = \int_{{\Bbb T}^2}^{\oplus} Ind_{\G_0}^{\G} \chi d\chi$$
Since $W^*(\G_0)$ is a masa of $W^*(\Gamma)$, it follows that
the induced representation $Ind_{\G_0}^{\G} \chi$ is irreducible
for almost all characters $\chi \in {\widehat\G_0}$
\cite [\S 8.5]{dixmier}.
}
\end{remark}
\bigskip 

\section{The Puk\'anszky invariant}

Let $A$ be a masa of a type $II_1$ factor $M$. Suppose that $M$ acts 
in its standard representation on the Hilbert space $L^2(M)$, which
is the completion of $M$ relative to the inner product defined by the
trace on $M$. Denote by $L$ the left regular representation of $M$ on 
$L^2(M)$ defined by $L_xy = xy$, and by $R$ the anti-representation defined
by $R_xy = yx$. Let $\mathcal A$ be the (abelian) von Neumann subalgebra
of ${\mathcal B}(L^2(M)$ generated by $L_A$ and $R_A$.
 A result of Popa \cite [Corollary 3.2]{po3} asserts that if 
$A$ is a regular masa then $\mathcal A$ is a masa of ${\mathcal B}(L^2(M))$.
Let $p_1$ denote the orthogonal projection of $L^2(M)$ onto the closed
subspace
generated by $A$. Then by \cite [Lemma 3.1]{po3},  $p_1$ is in the centre
of ${\mathcal A}'$ 
 and ${\mathcal A}'p_1$ is abelian. 
The Puk\'anszky invariant is the type of the (type $I$) von Neumann algebra
${\mathcal A}'(1 - p_1)$. It is an isomorphism invariant of the
pair $(A,M)$, since any automorphism of $M$ is implemented by a unitary
in ${\mathcal B}(L^2(M))$. The Puk\'anszky invariant has been computed for
some particular examples in \cite{puk}, \cite{rad}, \cite{br}.
 Popa showed that if ${\mathcal A}'(1 - p_1)$ is homogeneous of
type $I_n$ where $2 \le n \le \infty$ then $A$ is a singular masa
\cite [Remark 3.4]{po3}.

Now let $\Gamma$ be a group and $M = W^*(\Gamma)$. Then $L^2(M)$ is
naturally
identified with  $l^2(\Gamma)$ and $L$ and $R$ respectively become 
the left regular representation $L_g \delta_h = \delta_{gh}$
and the right regular anti-representation $R_g \delta_h = \delta_{hg}$ .

Let $\Gamma_0$ be an abelian subgroup of $\Gamma$ such that $A =
W^*(\Gamma_0)$ is a masa
of $M = W^*(\Gamma)$. Then ${\mathcal A} = (L_{\Gamma_0} \cup R_{\Gamma_0})''$.
 Denote by ${\mathcal D}$ the set of double cosets 
$D = \Gamma_0 g \Gamma_0 \in \Gamma_0 \setminus \Gamma / \Gamma_0$ which
have
the property that
\begin{equation} \label{3.1}
 g^{-1}\Gamma_0g \cap \Gamma_0 = \{ e \}
\end{equation}
Note that this condition depends only on the the double coset $\Gamma_0 g
\Gamma_0$
and not on the representative element $g$.

For each $D \in {\mathcal D}$ let $p_D$ denote the orthogonal projection from
$l^2(\Gamma)$
onto $l^2(D)$.  Since $l^2(D)$ is invariant under $L_{\Gamma_0}$ and
$R_{\Gamma_0}$
it follows that $p_D$ lies in the commutant ${\mathcal A}'$ of ${\mathcal A}$.

\begin{lemma} \label{3a}
If $C,D \in {\mathcal D}$ and $C \ne D$ then the projections $p_C$, $p_D$ are
 mutually orthogonal and equivalent in ${\mathcal A}'$.
\end{lemma}

{\sc Proof:} Let $C = \Gamma_0 c \Gamma_0$ and $D = \Gamma_0 d \Gamma_0$,
where
$c,d \in \Gamma$, and assume that $C \ne D$.
 Define a map $\phi : C \to D$ by $\phi (ucv) = udv$ for $u,v \in
\Gamma_0$.
The condition (\ref{3.1}) ensures that  $\phi$ is a well defined bijection.
 Moreover 
$$\phi(uk) = u\phi(k), \qquad \phi(kv) = \phi(k)v 
\qquad \quad \mbox{for}\quad u,v \in \Gamma_0, \quad  k \in C$$

Define a partial isometry $s$ on $l^2(\Gamma)$ by 
$$s\delta_k  =  \left\{ \begin{array}{ll}
                \delta_{\phi(k)} & \mbox{if} \quad k \in C,\\
                0  & \mbox{otherwise}
                \end{array}
\right.$$
Then $s^*s = p_C$ and $ss^* = p_D$.  We must show that $s \in {\mathcal A}'$.

If $k \in C$ and $u \in \Gamma_0$ then
$$L_us\delta_k = L_u\delta_{\phi(k)} = \delta_{u\phi(k)} =
\delta_{\phi(uk)}
= s\delta_{uk} = sL_u\delta_k$$

If $k \notin C$ and $u \in \Gamma_0$ then $uk \notin C$ and
$$L_us\delta_k = 0 = sL_u\delta_k$$
It follows that $L_us = sL_u$ for each $u \in \Gamma_0$.
Similarly $R_us = sR_u$ for each $u \in \Gamma_0$. Therefore $s \in {\mathcal
A}'$,
as required.

It is clear that $p_C$, $p_D$ are mutually orthogonal, since $C \cap D =
\emptyset$.
\qed

\begin{lemma} \label{3b} If $D = \Gamma_0 d \Gamma_0 \in {\mathcal D}$
 where $d \in \Gamma$, then $p_D$
is an abelian projection in ${\mathcal A}'$.
\end{lemma}

{\sc Proof:} (c.f. \cite[Lemma 4]{puk}.) By (\ref{3.1}) the map $(x,y)
\mapsto xdy$ defines
 a bijection
from $\Gamma_0 \times \Gamma_0$ onto $D$. Define a unitary operator 
$U : l^2(D) \to l^2(\Gamma_0 \times \Gamma_0)$ by
$(Uf)(x,y) = f(xdy)$. Then for $z_1,z_2 \in \Gamma_0$ and $f \in l^2(D)$,
$$(UL_{z_1}R_{z_2}f)(x,y) = f(z_1^{-1}xdyz_2^{-1}) =
(Uf)(z_1^{-1}x,yz_2^{-1}))$$

Thus $UL_{z_1}R_{z_2}U^{-1} = \Lambda_{z_1}P_{z_2}$ where the operators
on the right hand side are defined on $l^2(\Gamma_0 \times \Gamma_0)$ by
$(\Lambda_zf)(x,y) = f(z^{-1}x,y)$ and $(P_zf)(x,y) = f(x,z^{-1}y)$.
Therefore ${\mathcal A}p_D$ is spatially isomorphic to the subalgebra of 
${\mathcal B}(l^2(\Gamma_0 \times \Gamma_0))$ generated by the
operators $\Lambda_z$ and $P_z$, for $z \in \Gamma_0$.
Taking Fourier transforms, this algebra is itself spatially isomorphic to
the algebra $L^{\infty}(\widehat\Gamma_0 \times \widehat\Gamma_0)$ acting
by multiplication
on $L^2(\widehat\Gamma_0 \times \widehat\Gamma_0)$. The latter algebra is
maximal abelian
\cite [Vol.I, Example 5.1.6]{kr} and hence ${\mathcal A}p_D$ is maximal abelian
on $l^2(D)$.
It follows that $p_D{\mathcal A}'p_D$ is contained in ${\mathcal A}p_D$ and so is
abelian.
In other words $p_D$ is an abelian projection in ${\mathcal A}'$.
\qed

\begin{lemma} \label{3c}
Let $\Gamma$ be an $\widetilde A_2$ group and $\Gamma_0$ the abelian
subgroup
generated by three distinct commuting elements $a,b,c \in P$.
 Then the set $\mathcal D$ is infinite. That is,
there exists a sequence $\{ g_n \}$ of elements of $\Gamma$ such that

(i) $g_n^{-1}\Gamma_0g_n \cap \Gamma_0 = \{ e \}$ for $n = 1,2,\dots$

(ii) $\Gamma_0g_r\Gamma_0 \cap \Gamma_0g_s\Gamma_0 = \emptyset$ for $r \ne
s$
 
\end{lemma}

{\sc Proof:}
By Lemma \ref{a,b,c} we have
 $\lambda(a) \cap \lambda(b) = \{c\}$, $\lambda(b) \cap \lambda(c) = \{a\}$
and  $\lambda(c) \cap \lambda(a) = \{b\}$.
 Hence
$$ \# (\lambda(a) \cup \lambda(b) \cup \lambda(c)) = (q + 1) + q + (q - 1)
= 3q$$
Since the total number of points in $P$ is $q^2 + q + 1$ and
$q^2 + q + 1 > 3q$ for $q \ge 2$, we can choose  a generator $z_1 \in P$
with $z_1  \notin \lambda(a) \cup \lambda(b) \cup \lambda(c)$. Of course
this
implies in particular that $z_1 \notin \{a,b,c\}$.

Now choose $z_2 \in P$ such that $z_2 \notin \{ a,b,c \} \cup \lambda
(z_1)$.
This is possible since $\# (\{ a,b,c \} \cup \lambda (z_1)) \le q + 4 < 
q^2 + q + 1 = \#(P)$.  By induction we can choose a sequence
$z_2,z_3,\dots$
in $P$ such that $z_{j+1} \notin \{ a,b,c \} \cup \lambda (z_j)$.
Define $g_n = z_1z_2 \dots z_n$. 

Note that any element $g \in \Gamma_0$ has left normal form
$g = x^{-k}y^l$ where $x,y \in \{a,b,c\}$, $k,l \ge 0$ and $x \ne y$.

Suppose that (i) is not true. Then we can find an element 
$$g_n^{-1}x_1^{-m_1}y_1^{n_1}g_n = x_2^{-m_2}y_2^{n_2} \in
g_n^{-1}\Gamma_0g_n \cap \Gamma_0$$
where $x_i,y_i \in \{a,b,c\}$, $x_i \ne y_i$, $m_i,n_i \ge 0$ and $m_i +
n_i > 0$,  $i = 1,2$.
This means that
$$z_n^{-1}\dots z_1^{-1}x_1^{-m_1}y_1^{n_1}z_1\dots z_n =
x_2^{-m_2}y_2^{n_2}$$
where on each side of the equation at most one of the terms $x_i$,$y_i$  is
absent.
With or without such an absence, each side of the equation
is in left normal form, which contradicts uniqueness of left normal form.
This proves (i).

In order to prove (ii), it is enough to show that if
$g_r\in\Gamma_0g_s\Gamma_0$, then $r=s$.  Suppose that
$g_r = x_1^{-m_1}y_1^{n_1}g_sx_2^{-m_2}y_2^{n_2}$, where
$x_1^{-m_1}y_1^{n_1}$, $x_2^{-m_2}y_2^{n_2}$ are elements of
$\Gamma_0$ expressed as usual in left normal form. Then
$$z_1 \dots z_r = x_1^{-m_1}y_1^{n_1}z_1 \dots z_sx_2^{-m_2}y_2^{n_2}$$
This implies
$$x_1^{m_1}z_1 \dots z_ry_2^{-n_2} = y_1^{n_1}z_1 \dots z_sx_2^{-m_2}$$
Both sides are in right normal form, which
is unique \cite [I, Proposition 3.2]{cmsz}.
Since $z_j\notin\{a,b,c\}$ for any $j$, this implies that $r=s$.
\qed 

\begin{corollary} \label{3d}
Let $\Gamma$ be an $\widetilde A_2$ group and $\Gamma_0$ the abelian
subgroup
generated by three distinct commuting elements $a,b,c \in P$.
(So that $W^*(\Gamma_0)$ is a masa of $W^*(\Gamma)$, by Lemma \ref{z2}
and Corollary \ref{max}.)
Then the Puk\'anszky invariant contains a type $I_{\infty}$ summand.
\end{corollary}

{\sc Proof:} Lemmas \ref{3a}, \ref{3b} and \ref{3c} show that 
${\mathcal A}'(1 - p_1)$ contains an infinite family of mutually
orthogonal equivalent abelian projections.
\qed

\begin{remark} \label{later}
{\em
 Another consequence of Lemma \ref{3c}(i) is that 
 $W^*(\Gamma_0)$ has no nontrivial central sequences in $W^*(\Gamma)$. This
follows from \cite[Remark 4.2(2)]{po1}.
}
\end{remark}

The following result clarifies the situation in
a general group $\Gamma$.

\begin{proposition} \label{3e}
Let $\Gamma$ be an arbitrary group and $\Gamma_0$ an abelian subgroup
such that  $A = W^*(\Gamma_0)$ is a masa of $M = W^*(\Gamma)$.
Suppose that  $g^{-1}\Gamma_0g \cap \Gamma_0 = \{ e \}$ for all
$g \in \Gamma \setminus \Gamma_0$. Let $n = \#(\Gamma_0 \setminus \Gamma /
\Gamma_0-\{\Gamma_0\})$.
Then ${\mathcal A}'(1 - p_1)$ is homogeneous of type $I_n$.
\end{proposition}

{\sc Proof:} By assumption we have ${\mathcal D} = \Gamma_0 \setminus \Gamma /
\Gamma_0-\{\Gamma_0\}$.
 Lemmas \ref{3a} and \ref{3b} imply that $\{ p_D : D \in {\mathcal D} \}$  is a
family of
$n$ (possibly $n = \infty$) mutually orthogonal equivalent abelian
projections in
${\mathcal A}'$ with  sum $1 - p_1$. This proves the result.
\qed
\bigskip

Proposition \ref{3e} does not apply when $\Gamma$ is an $\widetilde A_2$
group.
However it does describe the  Puk\'anszky invariant for the abelian
subalgebra
generated by a generator of a free group as in Example \ref{0}.
\bigskip
\begin{corollary} \label{3f}
Let $\Gamma = \Gamma_0 \star \Gamma_1$ where $\Gamma_0 = \langle a \rangle
\cong {\Bbb Z}$
and $\Gamma_1$ is any nontrivial group. Then ${\mathcal A}'(1 - p_1)$ is
homogeneous of type 
$I_{\infty}$.
\end{corollary}

{\sc Proof:} Choose $b \in \Gamma_1 \setminus \{ e \}$.  The double cosets
$\Gamma_0 b \Gamma_0, \Gamma_0 bab \Gamma_0, \Gamma_0 babab \Gamma_0,
\dots$
are pairwise disjoint. Therefore $\#(\Gamma_0 \setminus \Gamma / \Gamma_0)
= \infty$. 
Clearly $g^{-1}\Gamma_0g \cap \Gamma_0 = \{ e \}$ for all
$g \in \Gamma \setminus \Gamma_0$.
The result therefore follows from Proposition \ref{3e}.
\qed
\bigskip

\section{Line central elements}

In this section we determine how many generators of $\G$ can
commute with a given generator.

\begin{proposition} \label{3} A generator $a \in P$ commutes with at most
$q+1$ other
generators in $P \setminus \{a\}$. Moreover if  $a$ commutes with  $q+1$ 
generators $x_1,x_2,\dots,x_{q+1}$ in $P\setminus \{a\}$ then
\begin{itemize}
\item [(i)]  $a \notin \lambda(a)$
\item [(ii)] $\{x \in P : ax = xa \} = \{a\} \cup \lambda(a)$.
\end{itemize}
\end{proposition}

{\sc Proof:} If $a$ commutes with $x \in P \setminus \{a\}$ then $x \in
\lambda(a)$
by Lemma \ref{2}. The  assertion (i) is immediate, since the line
$\lambda(a)$
contains $q+1$ points.
If  $a$ commutes with  $q+1$ 
generators $x_1,x_2,\dots,x_{q+1} \in P\setminus \{a\}$ then 
$\lambda(a) = \{x_1,x_2,\dots,x_{q+1}\}$ and assertions (i) and (ii) follow
from this.
\qed
\bigskip

 Suppose that a generator $a$ commutes with $q + 1$ other generators.
The following terminology comes from the fact that the points on the
 line $\lambda(a)$ then correspond to the generators which commute with
$a$.
 
\begin{definition} \label{line central} If $P$ is the set of generators of
a triangle presentation then an element $a \in P$ is called {\it line
central} if it commutes with $q+1$ elements of $P\setminus \{a\}$.
\end{definition}

\begin{remark} \label{examples} 
{\em
  Examination of the tables
 given at the end of \cite {cmsz} suggest
that line central elements may be relatively rare. However they do exist.
Here are some
examples.

(i) In the presentation (B.3) of \cite{cmsz}, the generator $a_0$ is line central
with $\lambda(a_0) = \{a_1,a_2,a_4\}$. In this case $q = 2$ and there are
seven generators
$a_0$, \dots , $a_6$. The relations involving $a_0$ are given by the
following
triples :
$(a_0 a_1 a_1),(a_0 a_2 a_2),(a_0 a_4 a_4)$.
This group acts on the building of $PGL(3,{\Bbb Q}_2)$.

(ii) The presentation (4.1) of \cite{cmsz} has an unusually large number of
line
central generators, namely $a_2$,$a_7$ and $a_9$ . In this case $q = 3$.
 The line corresponding to $a_2$ is
$\lambda(a_2) = \{a_3,a_4,a_7,a_9\}$. There are thirteen generators $a_0$,
\dots , $a_{12}$
and the relations involving  $a_2$ are: $(a_2 a_4 a_7),(a_2 a_7 a_4),(a_2
a_3 a_9),(a_2 a_9 a_3)$.
This group acts on the building of  $PGL(3,{\Bbb Q}_3)$.

(iii) The presentation (5.1) of \cite{cmsz} has exactly one line central
generator $a_2$.

(iv) In the presentation (64.1) of \cite{cmsz}, the generator $a_{12}$ is
line central. Again
$q = 3$.
The building on which this group acts is not that of any linear group
\cite{cmsz}.
}
\end{remark} 

 It may be significant that the groups which act on the building of
$PGL(3,{\Bbb F}_q((X)))$
do not appear to have a line central generator.
\medskip
\begin{proposition} \label{4} If an element $a \in P$ is line central and
$x \in P$ then
$$x \in \lambda(a) \Longleftrightarrow  a \in \lambda(x)$$
\end{proposition}

{\sc Proof:} Suppose $x \in \lambda(a)$. By Corollary \ref{3}, $ax = xa$.
Also
$axy = e$ for some $y \in P$ and so $xay = e$ as well. It follows that $a
\in \lambda(x)$.

Conversely, suppose that $a \in \lambda(x)$. Then $xay = e$ for some $y \in
P$. Therefore
$ayx = e$. This implies that $y \in \lambda(a)$.  By Corollary \ref{3}, $ay
= ya$.
Then $yax = e$, and so $x \in \lambda(a)$ 
\qed

\bigskip

\end{document}